\documentclass[reqno]{amsart}
\usepackage{amsmath}
\usepackage{amssymb}
\newtheorem{thm}{Theorem}[section]
\newtheorem{lemma}{Lemma}[section]
\newtheorem{prop}{Proposition}[section]
\newtheorem{cor}{Corollary}[section]

\theoremstyle{definition}

\theoremstyle{remark}

\newtheorem{remark}{Remark}
\def \mcv {\mathcal {V}}

\def \intm {\stackrel{\circ}{M}}

\def \mrn {{\mathbb R}^n}

\def \mrp {{\mathbb R}^p}
\def \mr {{\mathbb R}}

\def \mca {{\mathcal A}}

\def \mcv {{\mathcal V}}

\def \mc {{\mathbb C}}
\def \mh {{\mathbb H}}
\def \mn {{\mathbb N}}

\def \alf {\alpha}
\def \eps {\epsilon}
\def \la {\lambda}

\def \del {\delta}

\newcommand{\Id}{\operatorname{Id}}

\newcommand{\T}{\operatorname{T}}
\newcommand{\Tr}{\operatorname{Tr}}

\def \p {\partial}

\numberwithin{equation}{section}
\title[Inverse scattering on conformally compact manifolds]{Inverse
Scattering on conformally compact manifolds}
\author[Marazzi]{Leonardo Marazzi}
\address{Department of Mathematics
Purdue University,
West Lafayette IN 47907, U.S.A.}
\email{lmarazzi@math.purdue.edu}
\begin{document}
\input epsf\

\begin{abstract}
We study inverse scattering for $\Delta_g+V$ on $(X,g)$ a conformally compact manifold with metric $g,$ with variable sectional curvature $-\alf^2(y)$ at the boundary and    
$V\in C^\infty(X)$ not vanishing at the boundary. We prove that the scattering matrix  at a fixed energies $(\lambda_1,$ $\lambda_2)$  in a suitable subset of $\mc$,
determines $\alf,$ and the Taylor series of both the potential and the metric at the boundary.
\end{abstract}

\maketitle
\section{\bf Introduction}\label{intro}

In this note we study inverse scattering on conformally compact manifolds with non-constant asymptotic sectional curvatures. We work with $\Delta$ the negative Laplacian. We prove that the scattering matrix of $\Delta_g+V$, $g$
conformally compact, $V\in \mc^\infty$, at fixed energies $(\lambda_1,$ $\lambda_2)$ in a
suitable subset of $\mc$, determines the curvature $\alf,$ the Taylor series of  the potential $V$ at
the boundary and the Taylor series of $g$ near the boundary. This is a generalization and improvement of some results of \cite{jsb0}. It is a generalization to the case where the potential $V$ does not vanish at infinity and the curvature at infinity  is not constant. It is an improvement in the sense if we know a priori that $\alf_1=\alf_2$ we determine both Taylor series from the scattering matrix at just one energy.

The potentials considered here can be thought of as symbols of order zero. In the Euclidean setting scattering for such potentials has been studied by Hassel-Melrose-Vasy  \cite{HMV1,HMV2}, Sait$\overline o$ \cite{S}, Herbst 
\cite{H}, and Agmon-Cruz-Herbst \cite{ACH} among others. In this paper we get a first result on inverse scattering in the setting of conformally compact manifold with variable curvature at infinity. There does not seem to be any inverse result for the symbols which are potentials of order zero in the  Euclidean setting.

Let $X$ be a $C^\infty$ compact manifold of dimension $n+1$
with
boundary $\p X.$ We recall that  $x$ is a boundary defining function of $\p X$  if  $x\geq 0,$ $\p X=\{x=0\},$ and $dx\neq 0$ on $\p X.$ We assume $X$ is equipped with a Riemannian
metric $g$ such that for any defining function $x$ of $\p X$ the metric  $x^2g=\widetilde{g}$ is a $C^\infty$ non degenerate Riemannian metric up to $\p X$. The manifold $(X,g)$ is called a conformally compact manifold.

It is shown in \cite{Ma1} that if $\nu$ is the unit normal with respect to $\widetilde{g}$, $-(\p x/\p \nu)^2(y)=-\alf^2(y)$ are the sectional curvatures at the boundary.  When $(\p x/\p \nu)^2(y)$
is constant, the manifold is called asymptotically hyperbolic. The scattering theory in the setting of variable sectional curvature $\alf$ at infinity has been studied by Borthwick \cite{Borth}.

Following the proof of Lemma 2.1 of  \cite{Grah1} one can show that once fixed $\tilde g|_{\p X}$ there exists  a unique $C^\infty$ defining function $x$ of $\p X$, in a collar neighborhood $[0,\eps) \times \p X$ of  $\p X$, such that $|dx|_{\tilde g}=\alf.$ In this case we can write
\begin{gather}
g=\frac{dx^2}{\alf^2(y)x^2}+\frac{h(x,y,dy)}{x^2},  \;\ (x,y) \in [0,\eps) \times \p X.\label{hmet}
\end{gather}

Mazzeo and Melrose \cite{MM} studied the resolvent of the Laplacian for asymptotically
hyperbolic manifolds. They proved that the resolvent has a meromorphic
continuation to $\mc \backslash \{(1/2)(n-\mn_0)\}.$
Guillarmou \cite{CG}
proved that in general the resolvent may have essential singularities at $\{ (1/2)(n-\mn_0)\}.$
The
generalization to a variable curvature
at the boundary $\alf(y)$ was carried out by
Borthwick in \cite{Borth}; he proved the existence of the Poisson operator,
and meromorphic continuation of the resolvent, and the
scattering matrix.

By the spectral theorem the resolvent
$$
R_\lambda=\left( \Delta_g+V(x,y)-\lambda^2-\frac{n^2}{4}\right)^{-1}
$$
is well defined for $\Im \lambda<<0,$ and one would like to understand if it can meromorphically continued to a larger region of the complex plane.

The proof of  Proposition 5.3 of \cite{Borth} can be modified to our case to obtain

\begin{thm}\label{new1} Let  $V \in C^\infty(\overline{X}),$ and let $\la \in \mc\setminus (\Omega\cup D),$ with $\Omega$ defined in
\eqref{defgamma} below, and $D$ the discrete set of $\lambda\in\mc$ such that the resolvent $R_\lambda=\left( \Delta_g+V(x,y)-\lambda^2-\frac{n^2}{4}\right)^{-1}$ has a pole.  Let $x$ be such that \eqref{hmet} is satisfied. Given $f\in C^\infty(\p X),$ there exists a unique $u\in C^\infty(X),$ such that
\begin{gather*}\begin{gathered}
\left( \Delta_g+V(x,y)-\lambda^2-\frac{n^2}{4}\right) u(x,y)=0;\\
u(x,y)=x^{n-\sigma}F(x,y)+x^\sigma G(x,y),
\end{gathered}\end{gather*}
where $F,G \in C^\infty(\overline{X}),$ $F=f$ at $\p X$, and
$\sigma=\frac{n}{2}+\sqrt{\left(\frac{n}{2}\right)^2
-\frac{1}{\alf^2} (V(0,y)-\lambda^2-\frac{n^2}{4})}.$
\end{thm}

We outline the proof of \ref{new1} in section \ref{sec1}.
The Poisson operator is the map

\begin{gather}
\begin{gathered}
E_\lambda: C^\infty(\p X) \longrightarrow C^\infty(X) \\
E_\lambda:f\mapsto u,
\end{gathered}\label{poisson}
\end{gather}
and the scattering matrix $S_\lambda$  is defined by

\begin{gather}\begin{gathered}\label{defsc}
S_\lambda: C^\infty(\p X) \longrightarrow C^\infty(\p X) \\
S_\lambda:f\mapsto G\mid_{\p X}.
 \end{gathered} \end{gather}

In \cite{jsb0}, Joshi and  S\'{a} Barreto deal with the asymptotically
hyperbolic case and  show that for a fixed  $\lambda \in \mc \backslash Q,$ where $Q$ is a discrete set, the scattering matrix $S(\lambda)$
determines the Taylor series of the metric $g$ or the potential $V$, with the assumption that the
potential
vanishes at the boundary.  We carry out
the natural
extention of this approach to the conformally compact case  and for
potentials not vanishing
at the boundary.

Let $P_1$ and $P_2$ be the operators

\begin{gather}\begin{gathered}\label{operators}
 P_1=\Delta_{g_1}+V_1(x,y)-\lambda_1^2-\frac{n^2}{4},\\
P_2=\Delta_{g_2}+V_2(x,y)-\lambda_2^2-\frac{n^2}{4},
\end{gathered}\end{gather}
and we fix a product structure in which

\begin{gather}\label{V&g}
  \begin{gathered}
    g_j=\frac{dx^2}{\alf_j^2(y)x^2}+\frac{h_j(x,y,dy)}{x^2}\qquad i=1,2.
  \end{gathered}
\end{gather}

We introduce some notation necessary to state our Theorem. We denote by 
\begin{gather}\label{scpls}
D_i=\{\lambda\in \mc: (R_\lambda)_i=\left(\Delta_{g_i}+V_i(x,y)-\lambda^2-\frac{n^2}{4}\right)^{-1} \mbox{ has a pole}, i=1,2\},
\end{gather}
we also denote by
\begin{equation}\label{defgammai}
\Omega_i=\Omega_i'\cup \left[ \min_{\p X}\{ V_i(0,y)\}-\alf_{iM}^2 \frac{n^2}{4}+\frac{n^2}{4},
 \max_{\p X}\{V_i(0,y)\}-\alf_{im}^2\frac{n^2}{4}+\frac{n^2}{4}\right],\qquad i=1,2 
\end{equation}
where $\alf_m=\min_{\p X} \alf,$ $\alf_M=\max_{\p X} \alf,$ and 
\begin{equation*}
\Omega_i'= \left\{ \lambda \in \mc:\exists y\in \p X :  \sigma_i(\lambda,y) \in
\frac{n-\mn_0}{2} \right\},
\end{equation*}
with
$$
\sigma_i(\lambda,y)=\frac{n}{2}+\sqrt{\left(\frac{n}{2}\right)^2
-\frac{1}{\alf_i^2(y)} (V_i(0,y)-\lambda^2-\frac{n^2}{4})}.
$$
We denote by $S_1$ and $S_2$ the scattering matrices associated to $P_1$ and
$P_2$ respectively.

Our main Theorem is 

\begin{thm}\label{main}
Let $g_1,g_2$ and $V_1, V_2$ be as in \eqref{V&g}, $p\in \p X,$ and $\sigma=max\{\sigma_1,\sigma_2\}.$  Assume that near $p,$ 
$(S_1)_{\lambda_1}\equiv (S_2)_{\lambda_1}$  and $(S_1)_{\lambda_2}\equiv (S_2)_{\lambda_2}$ $\mod$
$\Psi_{2\max_{\p X}\Re\sigma -n-k}(\p X),$ $k\geq 1,$  $\lambda_j \in \mc\setminus(\Omega_1\cup \Omega_2\cup D_1\cup D_2)$ for $j=1,2,$ $D_i$ defined in  \eqref{scpls}.  Then
$\alpha_1=\alpha_2,$ $V_1(0,y)=V_2(0,y),$
$\Omega_1=\Omega_2,$ and
there is a  discrete set $Q \subset \mc\setminus (\Omega_1\cup D_1\cup D_2)$ such that if  $\lambda_1 \in \mc \backslash(\Omega_1\cup D_1\cup D_2 \cup Q),$   then  $h_1-h_2=O(x^k)$ near $p$, and
$V_2-V_1=O(x^k)$ near $p.$
\end{thm}
 Theorem \ref{main} can be restated  invariantly as
\begin{thm}\label{two}  Let $g_1,g_2$ and $V_1, V_2$ be as in \eqref{V&g}, $p\in \p X,$ and $\sigma=max\{\sigma_1,\sigma_2\}.$  Assume that near $p,$ 
$(S_1)_{\lambda_1}\equiv (S_2)_{\lambda_1}$  and $(S_1)_{\lambda_2}\equiv (S_2)_{\lambda_2}$ mod
$\Psi_{2\max_{\p X}\Re\sigma -n-k}(\p X),$ $k\geq 1,$  $\lambda_j \in \mc\setminus(\Omega_1\cup\omega_2\cup\cup D_1\cup D_2)$ for $j=1,2,$ $D_i$ defined in  \eqref{scpls}. Then
$\alpha_1=\alpha_2,$ $V_1(0,y)=V_2(0,y),$
$\Omega_1=\Omega_2,$ and
there is a  discrete set $Q \subset \mc\setminus (\Omega_1\cup D_1\cup D_2)$ such that if $\lambda_1 \in \mc \backslash(\Omega_1\cup \Omega_2\cup D_1\cup D_2 \cup Q),$  then $V_2-V_1=O(x^k)$ near $p$ and
there exists a diffeomorphism $\phi$ of a neighborhood $U\subset X$ of $p,$ such that $\phi^*g_1-g_2=O(x^{k-2})$ in $U.$
\end{thm}

In Section \ref{sec1}, we recall the definitions of the spaces of polyhomogeneous distributions of
\cite{Borth} which are needed
to carry out the analysis for the conformally compact geometry with variable
curvature at infinity. The reason for the introduction of these spaces
comes from the appearance of an indicial root which will depend on the space
variable $y$, through the boundary curvature function $\alf(y)$ and
the potential $V(0,y)$. In Section \ref{sec2}, we prove our main Theorem.

The author would like to thank his advisor Ant\^{o}nio S\'{a} Barreto for his
guidance and help on this paper, and anonymous referees for many helpful comments.

\section{\bf The Poisson Operator and the Scattering Matrix}\label{sec1}

In this section we prove Theorem \ref{new1}.

\subsection{\bf Boundary asymptotics}
In this subsection we recall the spaces of functions used in \cite{Borth}.

Let $M$ be a smooth manifold with corners, as defined in \cite{melcor}, and
let $\rho=(\rho_1,...,\rho_p)$ be the defining functions for the finitely many
boundary faces
$Y_1,...,Y_P$ of $M$. Let $\mcv_b(M)$ be the set of smooth vector fields
tangent to the boundary. We recall the space of conormal distributions
\begin{equation}
\mca^m=\{u\in C^\infty(\intm):\mcv_b^ku\in\rho^mL^\infty(\intm),  \forall k\},
\end{equation}
where $m\in \mrp$ and $\rho^m=\rho_1^{m_1}\cdot\cdot\cdot\rho_p^{m_p}$. And call the set
\begin{equation}
\mca^{m-}=\bigcap_{m'<m}\mca^{m'}.
\end{equation}

With this space defined, we recall for $\beta \in \mc^\infty(M;\mrp)$ the space of
polyhomogeneous distributions

\begin{equation}\label{conor}
\mca_\beta(M)=\{u\in C^\infty(M):\left[\prod_{l=0}^p\prod_{k=0}^{m_l-1}
(T_j-k)^{k+1} \right](\rho^{-\beta}u)
\in \mca^n(M), \forall n <m, \forall m \},
\end{equation}
where $T_j=\rho_j\p_{\rho_j}$. 

Lastly we recall the space of truncated expansion 
\begin{equation}
\mca_{\beta|q}(M)=\prod_{l=1}^p (\rho_l \ln \rho_l)^{q_l}\cdot\mca_\beta(M)=
\rho^{\beta}\left[\prod_{l=0}^p\prod_{k=0}^{q_l-1} (T_j-k)^{k+1}\right]
\rho^{-\beta}\cdot\mca_\beta(M).
\end{equation}

We refer the reader to \cite{Borth} for a more detailed description of the later
spaces and for a proof of the last equality. An important lemma which was proven in \cite{Borth}, tells us that these spaces only depend on the restriction
  to the boundary of $\beta.$ For our case $\beta$ will be the indicial root $\sigma$
that will be discussed next; it appears in the asymptotic
expansion that leads to the definition of the scattering matrix \eqref{defsc}.

\begin{lemma}\cite{Borth}\label{indofx}
The space $\mca_\beta$ is independent of the choice of $T_j$ and depends on
$\beta$ only through the restrictions $\beta|_{Y_i}$.
\end{lemma}

\subsection{\bf The indicial operator}\label{indop}
We adapt the parametrix construction of \cite{Borth}. For $g$ as in
\eqref{hmet}, we consider the Schr\"odinger operator
\begin{equation}
\Delta_g+V(x,y)-\lambda^2-\frac{n^2}{4}.
\end{equation}
We consider the indicial roots when restricted to the boundary $x=0.$ The indicial root $\sigma$ satisfy the equation
\begin{gather}\label{sigma}\begin{gathered}
-\alf^2\sigma(n-\sigma)+V(0,y)-\lambda^2-\frac{n^2}{4}=0\\
\Rightarrow \sigma_{\pm}=\frac{n}{2}\pm\sqrt{\left(\frac{n}{2}\right)^2
-\frac{1}{\alf^2} (V(0,y)-\lambda^2-\frac{n^2}{4})}.
\end{gathered}\end{gather}

We denote $\sigma=\sigma_+$, and therefore $\sigma_-=n-\sigma.$  Observe that $\sigma$ is holomorphic in $\lambda$ when

\begin{equation*}
\lambda^2 \notin \left[ \min_{\p X}\{ V(0,y)\}-\alf_M^2 \frac{n^2}{4}+\frac{n^2}{4},
 \max_{\p X}\{V(0,y)\}-\alf_m^2\frac{n^2}{4}+\frac{n^2}{4}\right],
\end{equation*}
where $\alf_m$ and $\alf_M$ are the minimum and maximum of $\alf$ at $\p X$ respectively.
Let
\begin{equation*}
\Omega'= \left\{ \lambda \in \mc:\exists y\in \p X :  \sigma(\lambda,y) \in
\frac{n-\mn_0}{2} \right\},
\end{equation*}
and then let
\begin{equation}
\Omega=\Omega'\cup \left[ \min_{\p X}\{ V(0,y)\}-\alf_M^2 \frac{n^2}{4}+\frac{n^2}{4},
 \max_{\p X}\{V(0,y)\}-\alf_m^2\frac{n^2}{4}+\frac{n^2}{4}\right], \label{defgamma}
\end{equation}
we have, just as in \cite{Borth} Lemma 3.2:

\begin{lemma}
Let $\lambda \in \mc \backslash (\Omega\cup D)$ then for $v \in \mca_{\sigma|1}$,
we can find
$u\in\mca_{\sigma|1}$ such that
\begin{equation*}
v-[\Delta_g+V(x,y)-\lambda^2-\frac{n^2}{4}]u \in \dot C^\infty(X).
\end{equation*}
\end{lemma}

This is the first ingredient of the parametrix construction in \cite{MM}. The
following corollary follows from the same arguments in \cite{Borth},

\begin{cor}
Let $\lambda \in \mc \backslash (\Omega\cup D)$, then for $f\in C^{\infty}$ there exists
$u\in \mca_\sigma$ such that
 \begin{gather*}
          \begin{gathered}
      \left[\Delta_g+V(x,y)-\lambda^2-\frac{n^2}{4}\right]u\in \dot C^\infty(X);\\
       u(x,y)\sim x^\sigma f(y) $ near $  x=0.
     \end{gathered}
\end{gather*}
\end{cor}
\subsection{\bf Stretched product}

We also recall the construction of the stretched product, which is
the manifold (with corners) obtained after blowing
up the product $X\times X$ along $\p \Delta \iota$, where $\p \Delta \iota=
(\p X \times \p X)\cap \Delta \iota \cong \p X$, and $\Delta \iota$ is the set of fixed
points of the involution $I$ that exchanges the two projections,

\begin{gather*}
I(\pi_L(X\times X))=\pi_r(X\times X).
\end{gather*}

Where $\pi_L(X\times X)$ is the projection onto the first component $X\times \p X$, and
$\pi_r(X\times X)$ the projection onto the second component $\p X\times X.$

We use the usual notation for the stretched product $X\times_0 X$ and denote
the blow-down map by:

\begin{gather}\label{b}
b:X\times_0 X\rightarrow X\times X.
\end{gather}

The process of blowing-up just described, amounts to the introduction of singular
coordinates near the corner, they are given near left face, in local projective coordinates,  by (with $Y=y-y'$)
\begin{gather}\label{lface}
s=\frac{x}{x'},\quad z=\frac{Y}{x'},\quad x',\quad y',
\end{gather}
near the front face by
\begin{gather}\label{fface}
\rho=\frac{x}{|Y|},\quad \rho'=\frac{x'}{|Y|},\quad r=|Y|,
\quad \omega=\frac{Y}{|Y|},\quad y,
\end{gather}
near the right face  by
\begin{gather}\label{rface}
t=\frac{x'}{x},\quad z'=-\frac{Y}{x},\quad x,\quad y.
\end{gather}

Setting
\begin{gather*}
R=\sqrt{(x')^2+x^2+|y-y'|^2}
\end{gather*}
the left, right, and front faces are characterized by $\rho=0$, $\rho'=0$, and
$R=0$ respectively.

\subsection{\bf Pseudodifferential operators}

We recall the class of pseudodifferential operators that we need. We are going to work on the space of half densities of the form

\begin{gather*}
\left|\frac{h(x,y)}{\alf(y)}\right|^{1/2}\left|\frac{dx}{x}
\frac{dy}{x^n}\right|^{1/2}, \qquad h\in \mc^\infty(X),\,
h \neq 0,\quad \alf \in \mc^\infty(\p X), \alf\neq 0.
\end{gather*}

The $C^\infty$ multiples of such a density are sections 
of the bundle of singular half densities $\Gamma^{1/2}_0(X).$ Similarly, and we refer to \cite{MM, jsb0} for the details, we can define the bundles  $\Gamma^{1/2}_0(X\times X),$ and $\Gamma^{1/2}_0(X\times_0 X).$

We can now recall the definition of the class of pseudodifferential operators
$^0\Psi^m(X,\Gamma_0^{1/2})$,  as
the aperators B
\begin{gather*}
B:\dot C^\infty(X;\Gamma_0^{1/2})\rightarrow C^{-\infty}(X;\Gamma_0^{1/2}),
\end{gather*}
having a Schwartz kernel
\begin{gather*}
K_B\in C^{-\infty}(X\times X;\Gamma_0^{1/2}),
\end{gather*}
whose lift to $X\times_0 X$ has a conormal singularity of order m at the lifted diagonal.

As in \cite{Borth}, define also  $\,^0\Psi_{\sigma_l,\sigma_r}(X\times_0 X,\Gamma_0^{1/2})$ to be
the class of operator whose (Schwartz) kernel satisfies
\begin{gather*}
b^*K\in \mca_{\sigma_l,\sigma_r,0}(X\times_0 X,\Gamma_0^{1/2}),
\end{gather*}
and are extendible across the front face. The residual class of the
construction is $\Psi_{\sigma_l,\sigma_r}$ the operator with kernels in
$\mca_{\sigma_l,\sigma_r}(X\times X,\Gamma_0^{1/2})$.

\subsection{\bf The resolvent, the Poisson operator and the scattering matrix}
We can now apply Proposition 4.2 of \cite{Borth} to use the parametrix
construction of \cite{MM} section 7  to get 
\begin{prop}\label{parametrix}
Let $\lambda \in \mc \backslash (\Omega\cup D)$, then there exists $M_\lambda$ and $F_\lambda$
holomorphic, such that
\begin{equation*}
[\Delta_g+V(x,y)-\lambda^2-\frac{n^2}{4}]M_\lambda=I-F_\lambda
 \end{equation*}
with $M_\lambda\in\, ^0\Psi^{-2}+\,^0\Psi_{\sigma_l,\sigma_r}$ and
$F_\lambda \in \Psi_{\infty,\sigma_r}$ is a compact operator.
\end{prop}

To apply analytic Fredholm theory we need the invertibility of $I-F_\lambda$ for at least one value of $\lambda.$ To do that one modifies the parametrix $M_\lambda$ by adding a smoothing operator of finite rank which guaranties that $(I-F_\lambda)^{-1}$ exists for $\lambda$ such that $\Re \lambda=0,$ and $\Im \lambda<<0.$ For the details of this construction we refer the reader to the second paragraph in the proof of Theorem 7.1 on page 301 of \cite{MM}. 

We decompose the resolvent as the pull-back using the blow-down
map b (that is $^0\Psi^m,\,$ $^0\Psi_{\sigma_l,\sigma_r}$),
and its residual class $(\Psi_{\sigma_l,\sigma_r})$ and  state this as a Proposition,

\begin{prop}\label{resolvent}
The resolvent
\begin{equation*}
R_\lambda=\left[\Delta_g+V(x,y)-\lambda^2-\frac{n^2}{4} \right]^{-1}:
\dot C^\infty(X)\rightarrow  C^\infty(\stackrel o X)
\end{equation*}
has a meromorphic continuation to $\lambda\in\mc \backslash(\Omega\cup D')$, and
structure
\begin{equation*}
R_\lambda\in\, ^0\Psi^{-2}+\,^0\Psi_{\sigma_l,\sigma_r}
+\Psi_{\sigma_l,\sigma_r}
\end{equation*}
\end{prop}
The proof of the existence of the Poisson operator and the scattering matrix
follow the
same as in \cite{Borth}, the Poisson operator is equal to
\begin{equation*}
E_\lambda=C(\lambda)x'^{-\sigma_r}R_\lambda \mid_{x'=0},
\end{equation*}

The following Proposition, which is proven in \cite{Borth}, is the final ingredient needed to prove Theorem \ref{new1},

\begin{prop}
For the Schwartz kernel of the Poisson operator
\begin{equation*}
E_\lambda f=\int_{\p X} E_\lambda(x,z')f(y')d\mu_{\p X} y'
\end{equation*}
and f$\in C^{\infty}(\p X)$, we have
\begin{equation*}
E_\lambda f\in \mca_{\sigma}(X)+\mca_{n-\sigma}(X)
\end{equation*}
\end{prop}

The proof of Theorem \ref{new1} follows; for the reader interested in the details we refer to \cite{Borth}.
The principal symbol of the scattering matrix is
\begin{equation*}
S_{\lambda}(\xi)=2^{n-2\sigma}\frac{\Gamma(n/2-\sigma)}{\Gamma(\sigma-n/2)}
|\xi|^{2\sigma-n}_{h_0}
\end{equation*}
for $\lambda \in \mc \backslash (\Omega\cup D')$.

\begin{remark} Notice that $S_\lambda$ is a pseudodifferential operator in $\Psi^m_{\eps,0}$ for every $\eps>0,$ and $m=2\max_{\p X}\Re\sigma-n$. 
\end{remark}

\section{\bf Proof of Theorem \ref{main}}\label{sec2}

We proceed to analyze the relationship
between scattering matrices 
associated to two distinct operators as in \eqref{operators}. Let's consider first the case where the scattering matrices agree at the
principal symbol level. In this case
\begin{multline}\label{scatmat}
\sigma_P (S_{1\lambda_i}(\xi))=2^{n-2\sigma_1(\lambda_i)}\frac{\Gamma(n/2-\sigma_1(\lambda_i))}{\Gamma(\sigma_1(\lambda_i)-n/2)}
|\xi|^{2\sigma_1(\lambda_i)-n}_{h_{10}}=\\
\sigma_P (S_{2\lambda_i}(\xi))=2^{n-2\sigma_2(\lambda_i)}\frac{\Gamma(n/2-\sigma_2(\lambda_i))}{\Gamma(\sigma_2(\lambda_i)-n/2)}
|\xi|^{2\sigma_2(\lambda_i)-n}_{h_{20}},
\end{multline}
where for $j=1,2$ and $i-1,2$
\begin{gather*}
\sigma_j(\lambda_i)=\frac{n}{2}+\sqrt{\left(\frac{n}{2}\right)^2
-\frac{1}{\alf_j^2(y)} (V_j(0,y)-\lambda_i^2-\frac{n^2}{4})}.
\end{gather*}
We use  that $|t\xi|_{hi_0}=t|\xi|_{hi_0}$ for every $t\in\mr$ to obtain

\begin{equation*}
2^{2\sigma_2-2\sigma_1}\frac{\frac{\Gamma(n/2-\sigma_1)}{\Gamma(\sigma_1-n/2)}
|\xi|^{2\sigma_1-n}_{h_{10}}}
{\frac{\Gamma(n/2-\sigma_2)}{\Gamma(\sigma_2-n/2)}|\xi|^{2\sigma_2-n}_{h_{20}}}
= t^{2(\sigma_1-\sigma_2)}.
\end{equation*}
This implies that $\sigma_1$ and  $\sigma_2$ are identical. Hence   
$|\xi|_{h_{20}}=|\xi|_{h_{10}}$ for every $\xi\neq 0$ and thus 
$h_{10}$ and $h_{20}$ are also equal.

Furthermore, using the equations for $\sigma_1$ and $\sigma_2$
at $\lambda_i$ $i=1,2:$ 
\begin{multline}
\sigma_1(\lambda_i)=\sigma_2(\lambda_i)\Rightarrow\\
\frac{n}{2}+\sqrt{\left(\frac{n}{2}\right)^2
-\frac{1}{\alf_1^2(y)} (V_1(0,y)-\lambda_i^2-\frac{n^2}{4})}=
\frac{n}{2}+\sqrt{\left(\frac{n}{2}\right)^2
-\frac{1}{\alf_2^2(y)} (V_2(0,y)-\lambda_i^2-\frac{n^2}{4})}.
\end{multline}
Which implies that 
\begin{gather*}
\frac{1}{\alf_1^2(y)} (V_1(0,y)-\lambda_i^2-\frac{n^2}{4})=
\frac{1}{\alf_2^2(y)} (V_2(0,y)-\lambda_i^2-\frac{n^2}{4}).
\end{gather*}
Rearranging 
\begin{gather*}
\frac{V_1(0,y)}{\alf_1^2(y)}-\frac{V_2(0,y)}{\alf_2^2(y)} =
(\lambda_i^2+\frac{n^2}{4})\left(\frac{1}{\alf_1^2(y)}-\frac{1}{\alf_2^2(y)}\right);
\end{gather*}
which implies $V_1(0,y)=V_2(0,y),$ and thus $\alf_1=\alf_2.$ The last to equalities imply that $\omega_1=\omega_2.$   

 Going back to the metrics, we have obtained
$h_1\mid_{\p X}=h_2\mid_{\p X},$ which means that there exists a tensor $L(y,dy)$ such that
\begin{gather*}
h_2=h_1+xL+O(x^2).
\end{gather*}
Next we obtain the higher order Taylor coefficients of the metric and potential from the lower order symbols of the scattering matrix.
We denote by 
$$
\delta_i=\det g_i, \mbox{ for } i=1,2.
$$
Just as in \cite{jsb0} we have
\begin{gather*}
\delta_2^{\pm 1/4}=\delta_1^{\pm 1/4}(1+x\cdot\frac{1}{4}\Tr(h_1(0,y)^{-1}L(0,y))+O(x^2)).
\end{gather*}

For the rest of the proof we only need to use that the scattering matrices agree at one energy, so we take $\lambda=\lambda_1$ and drop the subindex. Also 
\begin{gather*}
\sigma=\sigma(\lambda)=\frac{n}{2}+\sqrt{\left(\frac{n}{2}\right)^2
-\frac{1}{\alf^2(y_c)} (V(0,y_c)-\lambda^2-\frac{n^2}{4})}.
\end{gather*}
for the rest of the section. The fixed point $y_c$ will appear naturally after applying the normal operator.

\begin{remark}
Notice that if we assume that $\alf_1=\alf_2,$ one only need one energy to obtain $V_1(0,y)=V_2(0,y)$ from the previous argument.
\end{remark}

 Now assume $S_{1\lambda}-S_{2\lambda}\in\Psi_{2\max \Re \sigma-n-2}(\p X),$ and we want to show that $h_1(x,y,dy)-h_2(x,y,dy)=O(x^2),$ and $V_1-V_2=O(x^2).$ To do so we go further and get
information on the derivatives of $V$ and the metric $h$. Let $P_1$ and $P_2$ be as defined in \eqref{operators}. First we  compute
$P_2-P_1$ as in \cite{jsb0}\footnote{There is
a little correction to the computation in \cite{jsb0}, pointed out in
\cite{GSa}.}. The difference in the metric between our case and that of \cite{jsb0} is that $g_{00}=\frac{1}{\alf^2x^2}$, and
$\delta_i=\det|g_i|=\frac{det|h_1|}{(\alf_i(y)x^{n+1})^2}.$
 So the only term that will change module higher order terms in the computation of 
\begin{equation*}
\del^\frac{1}{4}\Delta_g(\del^{-\frac{1}{4}}f)=
\sum_{i,j=0}^n \del^{-\frac{1}{4}} \p_{z_i}(g^{ij}(f(\p_{z_j}\del^\frac{1}{4})-
\del^\frac{1}{4}(\p_{z_j}f))),
\end{equation*}
is the $i=j=0$ term. This term is equal to
\begin{equation*}
-x \alf_1^2\frac{(1-n)}{4}\T +
(\alf_2^2-\alf_1^2)\left( 2xf \p_x \ln \del^\frac{1}{4}+
x^2f(\p_x^2 \ln \del^\frac{1}{4}+ (\p_x \ln \del^\frac{1}{4} )^2)
+x^2\p_x^2 f + 2x\p_x f\right),
\end{equation*}
where $\T=\Tr(h_1(0,y)^{-1}L(0,y))$, and $\delta=\frac{\delta_2}{\delta_1}$.

Since $\alf_1=\alf_2$ and $V_1(0,y)-V_2(0,y)$ we have
\begin{equation}\label{symbol1}
P_2-P_1=x\left( \sum_{i,j=1}^{n}H_{ij}x\p _{y_i} x\p _{y_j}
-\alf_1^2\frac{(1-n)}{4}T \right) +
\sum_{j=1}^{\infty}x^j(V_2^{(j)}(0,y)-V_1^{(j)}(0,y))+x^{2}R.
\end{equation}
Where $H_{ij}=h_1^{-1}(0,y)L(0,y)h_1^{-1}(0,y),$
and $$
W^{(j)}=\left(\frac{\p^j}{\p x^j} V_2\right)(0,y)-\left(\frac{\p^j}{\p x^j}V_1\right)(0,y).
$$

To find the expansion on the difference of the scattering matrices we can proceed
as in \cite{jsb0}. Let $R_1$ and $R_2$ be the resolvents of $P_1$ and $P_2$
respectively, then

\begin{gather*}
P_2(R_1-R_2)=(P_2-P_1)R_1=x E R_1,
\end{gather*}
where $E$ is the right hand side of \eqref{symbol1} after factoring out an $x.$
To obtain $R_2$ as a perturbation of $R_1$ we need to find $F$ so that
\begin{gather*}
P_2(F)=x E R_1.
\end{gather*}
We set $x=x's$ and $F=x'F_1$, and as $x'$ commutes with $P_2$ we obtain
\begin{gather}\label{normal}
P_2(F_1)=s E R_1.
\end{gather}
At this step we use a fundamental tool developed in \cite{MM}, which is the normal operator. For the details of its construction and further properties we refer to \cite{MM} sections 2 and 5. We recall that the normal operator $N_P$ of $\Delta_g+V$ at a point $y_c\in\p X$ is given by $\alf^2(y_c)\Delta_0+V(0,y_c),$ where $\Delta_0$ is the Laplacian on $\mh^{n+1}.$ Here  we assume that the metric $\alf^2(y_c)h_0(y_c)$ is transformed by a linear change of variables to the identity. 

We apply the normal operator to \eqref{normal} to get
\begin{gather}\label{lambda}
\alf^2(y_c)(\Delta_0+V(0,y_c)-\lambda^2-\frac{n^2}{4})N_P F_1=N_P(s E R_1).
\end{gather}
The right hand side of \eqref{lambda} is in $\mca_{\sigma+1,\sigma-1}.$ We can now apply
Proposition 6.19 of \cite{MM} to deduce that $N_P F_1\in\mca_{\sigma,\sigma-1},$ thus by the mapping properties of $N_P,$ we can write $F=x'(F_1)=R\rho^\sigma \rho'^{\sigma}\gamma (\lambda);$ with $\gamma(\lambda)\in C^\infty(X\times_0X\backslash\Delta_0,\Gamma^{1/2}_0(X\times_0 X)).$ Next we follow
the construction of the expansion of $\gamma$ which applies just as in \cite{jsb0}: 

We recall Proposition 4.4 of \cite{jsb0}, which holds for our case and states that the kernel of $S_\lambda$ satisfies
$$
b_\p^*S_\lambda=\frac{1}{M_\sigma}
b^*(x^{-\sigma+n/2}((x')^{-\sigma+n/2})R_\lambda)|_{\rho=\rho'=0}.
$$
where   $b_\p$ is the blow-up of the manifold $\p X\times \p X$ along the
diagonal $\Delta\in \p X\times \p X$ (we refer the reader to \cite{jsb0} for
the details of this blow-up). Thus we can write 
$$
b_\p^*(S_{2\lambda}-S_{1\lambda})=\frac{1}{M_\sigma}
b^*(x^{-\sigma+n/2}(x')^{-\sigma+n/2}R\rho^\sigma \rho'^{\sigma}\gamma (\lambda))|_{\rho=\rho'=0}.
$$
In the coordinates $r=|y-y'|,$ $w=(y-y')/r,$ $y'$ 
$$ 
b_\p^*(S_{2\lambda}-S_{1\lambda})=r^{1-2\sigma+n}\gamma(\sigma,r,\omega,y,0,0)
\left|\frac{dr}{r}\frac{d\omega}{r^n}dy'\right|^{1/2}. 
$$
Taking the $r^{-n}$ factor out of the half-density we get that 
$\gamma(\sigma,r,\omega,y,0,0)
|d\omega dy'|$ is the restriction of 
$$
R^{n/2}\rho^{n/2}(\rho')^{n/2}\gamma(\lambda)\left|\frac{d\rho}{(\rho)^{n+1}}\frac{d\rho'}{(\rho')^{n+1}}\frac{dR}{R^{n+1}}d\omega dy'\right|^{1/2}
$$
to the intersection of the left, right and front face, which is $\rho=\rho'=R=0.$ Next we explicitly calculate what this is.
 
Using the blow-up coordinates
\begin{gather*}
s=x/x',\,\,\,\,\, z=(y-y')/x',
\end{gather*}
and that $V(0,y_c)-\lambda^2-\frac{n^2}{4}=\alf(y_c)^2\sigma(y_c)(n-\sigma(y_c)),$ $\Delta_0=-(x\p_x)^2+nx\p_x-(x)^2\Delta_{h_0(y)},$ the equation \eqref{lambda} transforms into
\begin{equation}
\alf(y_c)^2(-(s\p_s)^2+ns\p_s-(s)^2\Delta_{\Id (z)}+\sigma(y_c)(n-\sigma(y_c)))
(s^\sigma(1+s^2+|z|^2)^{\frac{1-2\sigma}{2}}\gamma (s,z))
=N_P(s E )G.
\end{equation} 
 
We drop the dependence of $y_c$ to simplify notation. It is well known (e.g. Lemma 2.1. \cite{GZ1}) that for the Green kernel of the operator $\Delta_0-\sigma(\sigma-n),$ acting on half-densities is given by
\begin{equation*}
G(s,z)=\left( \frac{\pi^{-\frac{n}{2}}}{2}
\frac{\Gamma(\sigma)}{\Gamma(\sigma-\frac{n-2}{2})}
\frac{s^\sigma}{(1+s^2+|z|^2)^\sigma}\right)
\left| \frac{ds}{s}\frac{dz}{s^n}dy'\right|^{1/2}+E_1.
\end{equation*}

Where $E_1$ has conormal singularity at ${s=1, z=0}.$ 
This can be used to compute the leading singularity of the kernel of
$S_2(\lambda) - S_1(\lambda)$ by computing the expansion for the restriction of 
$$
R^{n/2}\rho^{n/2}(\rho')^{n/2}\gamma(\lambda)\left|\frac{d\rho}{(\rho)^{n+1}}\frac{d\rho'}{(\rho')^{n+1}}\frac{dR}{R^{n+1}}d\omega dy'\right|^{1/2}
$$
to  $\rho=\rho'=R=0.$

Since $G$ acts as a convolution operator with respect to the group action defined in \cite{MM}, we have modulo $\left|\frac{ds}{s}\frac{dz}{s^n}dy'\right|^{1/2}$
\begin{equation}\label{maineq}
\alf^2(s^\sigma(1+s^2+|z|^2)^{\frac{1-2\sigma}{2}}\gamma (s,z))
=c(\sigma)\left(\sum_{i,j=1}^n H_{ij}(y_c)\p_{z_i}\p_{z_j}
I_1+ (W^{(1)}(y_c)-\alf_1\frac{1-n}{4}T(y_c))I_2\right) +\beta,
\end{equation}
where $\beta\in\mca_{\sigma,\sigma},$ and for $l=1,2$
\begin{equation}
I_l=I_l(\sigma,s,z)
=\int_0^\infty\int_{\mrn}
\frac{t^\sigma}{(1+t^2+|U|^2)^\sigma(1+s^2/t^2+|z-(s/t)U|^2)^\sigma}\left(\frac{s}{t}\right)^{\sigma+5-2l}\frac{dt}{t}dU.
\end{equation}
We rename $\alf^2(y_c)\gamma(s,z)$  as $\gamma(s,z),$ since for all the following computations $\alf^2(y_c)$ is just a constant and is not necessary for any purpose in the rest of the proof. 

We recall we are looking for the restriction of   
\begin{multline}
R^{n/2}\rho^{n/2}(\rho')^{n/2}\gamma(\lambda)\left|\frac{d\rho}{(\rho)^{n+1}}\frac{d\rho'}{(\rho')^{n+1}}\frac{dR}{R^{n+1}}d\omega dy'\right|^{1/2}=\\
(x')^{n/2}\frac{s^{n/2}}{(1+s+|z|^2)^{n/4}}\gamma(\lambda)\left|\frac{dx'}{(x')^{n+1}}\frac{ds}{s}\frac{dz}{s^n}dy'\right|^{1/2},
\end{multline}
to $\{x'=s=0, |z|=\infty\ \}.$ First we restrict to $\{ x'=0\}$ to get
$$
\frac{s^{n/2}}{(1+s+|z|^2)^{n/4}}\gamma(\lambda)\left|\frac{ds}{s}\frac{dz}{s^n}dy'\right|^{1/2}.
$$
Which when restricted to $\{s=0, |z|=\infty\}$ is the same as
\begin{equation}\label{restric}
\frac{s^{n/2}}{|z|^{n/2}}\gamma(\lambda)\left|\frac{ds}{s}\frac{dz}{s^n}dy'\right|^{1/2}
\end{equation}
restricted to $\{s=0, |z|=\infty\}.$
Using $\omega=z/|z|$ \eqref{restric} becomes $A(\omega)|d\omega dy'|^{1/2}.$

By \eqref{maineq} $A(\omega)$ is given by 
\begin{multline}
A(\omega)=\lim_{s\rightarrow 0, |z|\rightarrow\infty}\gamma(s,z)=\\
\lim_{s\rightarrow 0, |z|\rightarrow\infty}
\frac{1}{s^\sigma(1+s^2+|z|^2)^{\frac{1-2\sigma}{2}}}
\left[\sum_{i,j=1}^n H_{ij}(y_c)\p_{z_i}\p_{z_j}
I_1+ (W^{(1)}(y_c)-\alf_1\frac{1-n}{4}T(y_c))I_2\right].
\end{multline}  
Setting $|z|u=s/t$ and $U=(t/s)|z|\tilde V.$ Then as in Lemma 5.1 of \cite{jsb0}
\begin{equation}
I_l(\sigma,s,z)
=s^\sigma|z|^{-2\sigma+5-2l}\int_0^\infty\int_{\mrn}
\frac{t^{2\sigma+4-2l-n}}{(u^2+s^2/|z|^2+|V|^2)^\sigma(1/|z|^2+u^2+|\omega-\tilde V|^2)^\sigma}dVdu.
\end{equation}
We set
\begin{equation}
T_l(\sigma,s,z)
=\int_0^\infty\int_{\mrn}
\frac{t^{2\sigma+4-2l-n}}{(u^2+s^2/|z|^2+|V|^2)^\sigma(1/|z|^2+u^2+|\omega-\tilde V|^2)^\sigma}dVdu.
\end{equation}

We rotate $\omega$ to $e_1=(1,0,...,0).$ The next step is to analyze the limits 
$$
\lim_{s\rightarrow 0, |z|\rightarrow\infty}\frac{1}{s^\sigma(1+s^2+|z|^2)^{\frac{1-2\sigma}{2}}}I_2,
$$
and 
$$
\lim_{s\rightarrow 0, |z|\rightarrow\infty}\frac{1}{s^\sigma(1+s^2+|z|^2)^{\frac{1-2\sigma}{2}}}
\p_{z_i}\p_{z_j}I_1.
$$
To do that we  recall the lemma from \cite{jsb0}
\begin{lemma}\cite{jsb0}
For $k\geq 1,$ and for $2\Re \sigma \geq\max\{n-k+1,k+2\}.$ Let 
$$
J(l,k,\sigma)=\int_0^\infty\int_{\mrn}\frac{u^{2\Re\sigma+k+3-2l-n}}{(u^2+|v|^2)^{\Re\sigma}(u^2+(e_1-v)^2)^{\Re\sigma}}dvdu,
$$
where $l=1,2;$ $e_1=(1,0,...,0).$ Then $|J(l,k,\sigma)|<\infty.$
\end{lemma}

We can apply dominated convergence to get
$$
T_2(\sigma)=\lim_{s\rightarrow 0, |z|\rightarrow\infty} s^{-\sigma}|z|^{2\sigma-1}I_2(\sigma,s,z)
=
\int_0^\infty\int_{\mrn}\frac{u^{2\sigma-n}}{(u^2+|v|^2)^{\sigma}(u^2+(e_1-v)^2)^{\sigma}}dvdu.
$$
Noticing that
\begin{equation}
\p_{z_j}\p_{z_i}I_1(s,z)=\tilde C_1(\sigma)s^\sigma(\p_{z_j}\p_{z_i}|z|^{3-2\sigma})T_1(s,z)+O(s^{\Re \sigma}|z|^{-1-2\Re \sigma}).
\end{equation}
We have 
\begin{multline}\label{penultima}
\lim_{s\rightarrow 0, |z|\rightarrow\infty}\frac{1}{s^\sigma(1+s^2+|z|^2)^{\frac{1-2\sigma}{2}}}
\p_{z_j}\p_{z_i}I_1(s,z)=\\
\lim_{s\rightarrow 0, |z|\rightarrow\infty}\frac{1}{s^\sigma(1+s^2+|z|^2)^{\frac{1-2\sigma}{2}}}\tilde C_1(\sigma)s^\sigma(\p_{z_j}\p_{z_i}|z|^{3-2\sigma})T_1(s,z).
\end{multline}
Changing coordinates from $z=(y-y')/x$ to $Y=(y-y'),$ we substitute 
\begin{equation}
\p_{z_j}\p_{z_i}|z|^{3-2\sigma}=|z|^{1-2\sigma}|Y|^{2\sigma-1}\p_{Y_j}\p_{Y_i}|Y|^{3-2\sigma}
\end{equation}
into equation \eqref{penultima} to get 
\begin{multline}
\lim_{s\rightarrow 0, |z|\rightarrow\infty}\frac{1}{s^\sigma(1+s^2+|z|^2)^{\frac{1-2\sigma}{2}}}
\p_{z_j}\p_{z_i}I_1(s,z)=\\
\frac{1}{s^\sigma(1+s^2+|z|^2)^{\frac{1-2\sigma}{2}}} s^\sigma \tilde C_1(\sigma)|z|^{1-2\sigma}|Y|^{2\sigma-1}\p_{Y_i}\p_{Y_j}|Y|^{3-2\sigma}\lim_{s\rightarrow 0, |z|\rightarrow\infty}T_1(s,z)=\\
C_1(\sigma )|Y|^{2\sigma-1}\p_{Y_i}\p_{Y_j}|Y|^{3-2\sigma}T_1(\sigma) 
\end{multline}
With
$$
T_1(\sigma)=\int_0^\infty\int_{\mrn}\frac{u^{2\sigma-2-n}}{(u^2+|v|^2)^{\sigma}(u^2+(e_1-v)^2)^{\sigma}}dvdu.
$$

This way we obtain the form for the
leading singularity of
$S_2(\lambda)-S_1(\lambda),$ which is 

\begin{multline}\label{lastprinc}
\gamma(\sigma,0,\omega,y,0,0)=\gamma(\sigma,0,Y/|Y|,y,0,0)=\\
C(\sigma)
\left(T_1(\sigma)\sum_{i,j=1}^\infty H_{ij}(y_c)|Y|^{2\sigma-1}\p_{Y_i}\p_{Y_j}|Y|^{3-2\sigma}
+T_2(\sigma)(W^{(1)}(y_c)
-\alf_1^2(y)\frac{1}{4}\cdot(1-n)T(y_c))\right).
\end{multline}
times a non-vanishing smooth half-density. 
By assumption  \eqref{lastprinc} is equal to zero for every $Y\neq 0.$ We take $Y=(1,0,...,0)$ and get that 
\begin{equation}
H_{1,1}(y_c)=-\frac{T_2(\sigma)}{(1-\sigma)T_1(\sigma)}(W^{(1)}(y_c)
-\alf_1^2(y)\frac{1}{4}\cdot(1-n)T(y_c)).
\end{equation}
The same way taking $Y=(0,1,0,...,0)$ we get that $H_{1,1}(y_c)=H_{2,2}(y_c),$ and in the same way we get that  
$H_{i,i}(y_c)=H_{j,j}(y_c)$ for every $i,j=1,...,n.$ We now put 
$Y=(1,1,0,...,0),$  into \eqref{lastprinc} and set it equal to zero to get
\begin{multline}
H_{1,1}(y_c)((3-2\sigma)2^{3/2-\sigma}+(3-2\sigma)(1-2\sigma)2^{1/2-\sigma}))+H_{1,2}(3-2\sigma)(1-2\sigma)2^{1/2-\sigma})=\\-\frac{T_2(\sigma)}{T_1(\sigma)}(W^{(1)}(y_c)
-\alf_1^2(y)\frac{1}{4}\cdot(1-n)T(y_c)).
\end{multline}
The same way taking $Y=(1,0,1,...,0),$ we get $H_{1,2}=H_{1,3}.$ The process can be continued to get that $H_{i,j}=H_{l,k}$ for all $i,j,k,l=1,...,n.$
Setting equation \eqref{lastprinc} to be zero looks now like
\begin{equation}\label{last}
T_1(\sigma)H_{1,1}(y_c)\sum_{i,j=1}^\infty (y_c)|Y|^{2\sigma-1}\p_{Y_i}\p_{Y_j}|Y|^{3-2\sigma}
+T_2(\sigma)(W^{(1)}(y_c)
-\alf_1^2(y)\frac{1}{4}\cdot(1-n)T(y_c)=0.
\end{equation}
Thus  we take $\lambda\in \mc \backslash(\Omega_1\cup D_1\cup D_2 \cup Q)$ with $Q$ the union of the discrete set $D_1$ of zeros of $T_1(\sigma),$ with  the discrete set $D_2$ of zeros of $T_2(\sigma),$ and the discrete set $D_3$ of zeros of $C(\sigma).$ We have that \eqref{last} can only happen if $H(y_c)=0,$ i.e. $L(y_c)=0;$ which also implies that $V^{(1)}_2(0,y_c)=V^{(1)}_1(0,y_c).$

The same argument can be applied in exactly the same way to get that the $V_i$s and $H_i$s agree to higher order, obtaining Theorem \ref{main}.

\end{document}